\documentclass[11pt,reqno]{amsart}
\usepackage[paper=a4,DIV=12,BCOR=16mm]{typearea}
\usepackage[T1]{fontenc}
\usepackage[english]{babel}
\usepackage{graphicx,float}
\newcommand{\CC}{\mathbb{C}}
\newcommand{\EE}{\mathbb{E}}
\newcommand{\NN}{\mathbb{N}}
\newcommand{\RR}{\mathbb{R}}
\newcommand{\ZZ}{\mathbb{Z}}
\theoremstyle{plain}
\newtheorem{theorem}{Theorem}[section]
\newtheorem*{corollary*}{Corollary}
\newtheorem{lemma}[theorem]{Lemma}
\theoremstyle{remark}
\newtheorem{remark}[theorem]{Remark}
\begin{document}
\title[Scale-free uncertainty principles]{Scale-free uncertainty principles and Wegner estimates for random breather potentials \\[2ex] 
Principes d'incertitude ind\'ependants de l'\'echelle et estim\'ees de Wegner pour des potentiels ``random breather''
}
\author[I.~Naki\'c]{Ivica Naki\'c}
\address[I.N.]{University of Zagreb, Croatia}
%
\author[M.~T\"aufer]{Matthias T\"aufer}
\address[M.T., M.T. \& I.V.]{Chemnitz University of Technology, Germany}
%
\author[M.~Tautenhahn]{Martin Tautenhahn}
%
\author[I.~Veseli\'c]{Ivan Veseli\'c}

\thanks{
 \copyright 2014 by the authors. Faithful reproduction of this article is permitted for non-commercial purposes.
{\today, \jobname.tex}}

\thanks{
This work has been partially supported by the DFG under grant \emph{Eindeutige Fortsetzungsprinzipien und Gleichverteilungseigenschaften von Eigenfunktionen}.
Part of these interactions have been supported by the binational German-Croatian 
DAAD project \emph{Scale-uniform controllability of partial differential equations}. Moreover, I.N.\ was partially supported by HRZZ project grant 9345.} 

\begin{abstract}
We present new scale-free quantitative unique continuation principles for 
Schr\"odinger operators. They apply to linear combinations of eigenfunctions corresponding to eigenvalues  
below a prescribed energy, and can be formulated as an uncertainty principle for spectral projectors.
This extends recent results of Rojas-Molina \& Veseli\'c \cite{RojasMolinaV-13}, and Klein \cite{Klein-13}. We apply the
scale-free unique continuation principle to obtain a Wegner estimate for a random Schr\"odinger operator of breather type.
It holds for arbitrarily high energies.  
Schr\"odinger operators with random breather potentials have a non-linear dependence on random variables.
We explain the challenges arising from this non-linear dependence.

%
\mbox{}\\
\textsc{Resum\'e.} Nous pr\'{e}sentons de nouveaux principes de continuation unique ind\'ependants de l'\'echelle pour des op\'erateurs de Schr\"odinger.
Nos r\'esultats concernent des combinaisons lin\'eaires de fonctions propres correspondant aux valeurs propres au-dessous d'une \'energie prescrite et ils peuvent \^etre formul\'es en terme de principes d'incertitude pour des projecteurs spectraux.
Ceci g\'en\'eralise des r\'esultats  r\'ecents de Rojas-Molina \& Veseli\'c \cite{RojasMolinaV-13}, et Klein \cite{Klein-13}.
Nous utilisons des estimations de continuation unique ind\'ependantes de l'\'echelle et obtenons ainsi une estimation de Wegner pour un op\'erateur de Schr\"odinger al\'eatoire de type ``breather''.
De tels op\'erateurs d\'ependent des variables al\'eatoires d'une fa\c{c}on non-lin\'eaire et nous expliquons les difficult\'es li\'ees \`a cette non-lin\'earit\'e.
\end{abstract}
\maketitle
\section{Introduction}
A Wegner estimate is an upper bound on the expected number
of eigenvalues in a prescribed energy interval of a finite box Schr\"odinger operator.
The expectation here refers to the potential which is random. 
The most studied example in this situation is the so-called 
alloy-type potential, sometimes also called continuum Anderson model, cf.~Remark \ref{r-nonlinearity} below. A particular feature of this model is that 
randomness enters the model via a countable number of random variables, 
and these random variables influence the potential in a linear way.
In the random breather model we study here, this dependence is no longer linear but becomes non-linear.
What remains is the monotone dependence of the potential on the random variables.
The topic of the present note is to explain how to effectively use this monotonicity
in order to derive a Wegner estimate.
This only works if it is possible to cast the monotonicity in a quantitative form.
\par
In order to achieve this, we use a scale-free uncertainty relation or unique continuation principle for spectral projectors of Schr\"odinger operators, presented in Theorem~\ref{thm:sfuc}. A proof of Theorem~\ref{thm:sfuc} will be given in the forthcoming paper \cite{NakicTTV-prep}. It answers positively a question raised in \cite{RojasMolinaV-13}. A partial answer (for small energy intervals) had been given shortly after in \cite{Klein-13}.
Previously, there has been in the literature on random operators a plethora of related results, 
applicable in specialised situations, see e.g.\ \cite{RojasMolinaV-13} for a discussion.
However, the lack of a result like  Theorem~\ref{thm:sfuc} was a bottleneck for further progress.
\par
Estimates as in Theorem~\ref{thm:sfuc} have been developed and applied in a different area of mathematics, namely control theory for partial differential equations, 
starting with the seminal paper  \cite{LebeauR-95}. In this context they are called \emph{spectral inequalities}.
In fact, our proof of Theorem~\ref{thm:sfuc} highlights how ideas from the theory of random Schr\"odinger operators and control theory complement each other in an efficient way.
\section{Results}
Let $d \in \NN^\ast = \{1,2,3,\ldots\}$, $\delta>0$, $L\in \NN^\ast$ and $V\colon \RR^d \to \RR$ measurable and bounded.
Denote by $\Lambda_L = (-L/2 , L/2)^d$ a cube in $\RR^d$, by $$S_{L,\delta} = \Lambda_L \cap \Big(\bigcup\limits_{j \in \ZZ^d} B(z_j , \delta) \Big)$$ 
 the union of $\delta$-balls centered at  the points $z_j$ and contained in the corresponding unit cubes $\Lambda_1+j$, and by $H_L$ one of the self-adjoint restrictions of the Schr\"odinger operator 
$-\Delta + V$ to $\Lambda_L$ with either Dirichlet, Neumann, or periodic boundary conditions.
We formulate a scale-free quantitative unique continuation property for the operator $H_L$.
\begin{theorem}\label{thm:sfuc}
There is $K_0 = K_0 (d)$ such that for all $\delta\in (0,1/2)$, all $E \in \RR$, all measurable and bounded $V : \RR^d \to \RR$, all $L \in \NN^\ast$, all sequences $(z_j)_{j \in \ZZ^d} \subset \RR^d$ such that $\forall j \in \ZZ^d: B(z_j , \delta) \subset \Lambda_1 + j $ and all linear combinations of eigenfunctions
\[
\psi= \sum\limits_{n \in \NN^\ast \colon E_n\leq E} \alpha_n \psi_n 
\]
(where $\psi_n \in W^{2,2} (\Lambda_L ; \RR)$ form an orthonormal basis and satisfy $H_L \psi_k = E_n\psi_n$ and $(\alpha_n)_{n \in \NN^\ast}$ is a sequence in  
$\CC$) we have
\[
\int_{S_{L,\delta}} \lvert \psi \rvert^2  \geq C_{\rm sfuc} \int_{\Lambda_L} \lvert \psi \rvert^2 , 
\quad \text{where} \quad
C_{\rm  sfuc} = \delta^{K_0 \bigl(1 + \lVert V \rVert_\infty^{2/3}+ \lvert E \rvert^{1/2}\bigr)} .
\]
\end{theorem}
The constant $C_{\rm  sfuc}$ is called an observability constant or a scale-free unique continuation constant.
We can reformulate this statement as an uncertainty principle. For this purpose, denote by $\chi_I(H_L)$ the spectral projector of $H_L$ onto an interval $I$
and  by $W_{L,\delta}$ the characteristic function of the set $S_{L,\delta}$.
\begin{corollary*}
Under the same assumptions as in the above Theorem we have
\begin{equation}\label{eq:abstractUCP}
\chi_{(- \infty , E]}(H_L) \ W_{L,\delta} \ \chi_{(- \infty , E]}(H_L) 
\geq \delta^{K_0 (1 + \lVert V \rVert_\infty^{2/3}+ \lvert E \rvert^{1/2} )} \chi_{(- \infty , E]}(H_L) .
\end{equation}
\end{corollary*}
Inequality~\eqref{eq:abstractUCP} is to be understood in the sense of quadratic forms.
\begin{remark}
For our purposes, the explicit quantitative dependence of the constant $C_{\rm  sfuc} = C_{\rm  sfuc}(\delta, \lVert V \rVert_\infty, E)$ is essential.
In particular,  $C_{\rm sfuc}$ \emph{does not} depend on the scale $L \in  \NN^\ast$. 
It depends on the radius $\delta$ in a polynomial way, and on $\lVert V \rVert_\infty$ and $\lvert E \rvert$ in an exponential way. 
Note also that the constant is unaffected by a translation of a ball $B(z_j , \delta)$ as long as it stays in the corresponding unit cube.
For any $K_V \geq 0$, the bound is uniform in the ensemble of potentials $\{V\colon \RR^d \to [-K_V , K_V] \text{ measurable} \}$.
This is important, because we want to apply the theorem to random Schr\"odinger operators. 
There the constant must not depend on the particular configuration of randomness.
Since the operator is lower bounded, we have $\chi_{(- \infty , E]}(H_L) =\chi_{[- \lVert V \rVert_\infty, E]}(H_L)$.
\end{remark}
To put this into context let us cite similar results from \cite{LebeauR-95,RousseauL-12,BourgainK-13}. First we cite a special case of \cite[Corollary 2]{LebeauR-95} and \cite[Theorem 5.4]{RousseauL-12}.
\begin{theorem}[\cite{LebeauR-95}] \label{thm:LebeauR-95} Let $\Omega \subset \RR^d$ be bounded, open and connected, $T > 0$ and $\omega \subset \Omega \times [0,T]$ open and non-empty with $\overline{\omega} \subset (0,T) \times \mathring{\Omega}$. Then there is  $C = C (T,\Omega,\omega) > 0$ such that
\[
 \forall \psi \in L^2 (\Omega) : \quad
 \lVert \mathrm{e}^{T \Delta} \psi \rVert_{L^2 (\Omega)}^2 
 \leq C \iint_\omega \lvert \mathrm{e}^{t \Delta} \psi \rvert^2 .
\]
\end{theorem}
While this result applies to parabolic equations, the next one is an adaptation to the elliptic setting.
\begin{theorem}[\cite{RousseauL-12}]\label{thm:RousseauL-12}
 Let $\Omega \subset \RR^d$ be bounded, open and connected, and $\omega \subset \Omega$ open and non-empty with $\overline \omega \not = \Omega$. Then there is $K = K (\omega , \Omega) > 0$ such that for all sequences $(\alpha_j)_{j \in \NN^\ast} \subset \CC$ and all $E > 0$ we have
 \[
  \Bigl
  \lVert \sum\limits_{n \in \NN^\ast \colon E_n\leq E} \alpha_n \psi_n 
  \Bigr\rVert^2_{L^2 (\Omega)} 
  \leq K \mathrm{e}^{K \sqrt{E}}\,
  \Bigl\lVert\sum\limits_{n \in \NN^\ast \colon E_n\leq E} \alpha_n \psi_n \Bigr\rVert^2_{L^2 (\omega)} .
 \]
Here, $E_n$, $n \in \NN^\ast$, denote the ordered eigenvalues of $-\Delta$ on $\Omega$ with Dirichlet boundary conditions with corresponding eigenfunctions $\psi_n$, $n \in \NN^\ast$.
\end{theorem}
In contrast to Theorem~\ref{thm:sfuc} the dependence of the observability constant in Theorems~\ref{thm:LebeauR-95} and \ref{thm:RousseauL-12} on the geometry of $\Omega$ and $\omega$ is not known.
Next we cite \cite[Theorem 3.4]{BourgainK-13}, where a quantitative dependence on the observability constant similar to Theorem~\ref{thm:sfuc} is obtained. 
It applies to approximate solutions of the stationary Schr\"odinger equation.
A common feature of Theorems \ref{thm:BourgainK-13} and \ref{thm:sfuc} is the appearance of the term $K^{2/3}$ and $\lVert V \rVert_\infty^{2/3}$, respectively in the exponent. 
This is due to the use of Carleman estimates.
\par%
\begin{theorem}[\cite{BourgainK-13}] \label{thm:BourgainK-13}
 Let $\Omega \subset \RR^d$ be an open subset of $\RR^d$ and consider a real measurable function $V$ on $\Omega$ with $\lVert V \rVert_\infty \leq K < \infty$. 
 Let $\psi \in W^{2,2}(\Omega)$ be real-valued and $\xi \in L^2(\Omega)$ be defined by $-\Delta \psi + V \psi = \xi$ almost everywhere on $\Omega$. 
 Let $\Theta \subset \Omega$ be a bounded and measurable set where $\lVert \psi \rVert_{L^2(\Theta)} > 0$. Set
 \[
  \mathcal{Q}(x, \Theta) := \sup_{y \in \Theta} \lvert y - x \rvert \quad \text{for}\ x \in \Omega.
 \]
 Consider $x_0 \in \Omega \setminus \overline\Theta$ such that $\mathcal{Q} = \mathcal{Q}(x_0, \Theta) \geq 1$ and $B(x_0, 6 \mathcal{Q} + 2 ) \subset \Omega$. Then given $0 < \delta \leq \min \{ \operatorname{dist}(x_0, \Theta), 1/24 \}$, we have
 \[
  \left( \frac{\delta}{\mathcal{Q}} \right)^{m(1 + K^{2/3}) (\mathcal{Q}^{4/3} + \log \frac{ \lVert \psi \rVert_{L^2(\Omega)}}{\lVert \psi \rVert_{L^2(\Theta)}} ) } 
  \lVert \psi \rVert_{L^2(\Theta)}^2 \leq \lVert \psi \rVert_{L^2(B(x_0, \delta))}^2 + \delta^2 \lVert \xi \rVert_{L^2(\Omega)}^2 ,
 \]
where $m > 0$ is a constant depending only on $d$.
\end{theorem}
Now we discuss an application of Theorem~\ref{thm:sfuc} to \emph{random breather models},
a class of random Schr\"odinger operators where the randomness enters the potential in a non-linear way.
Consider a sequence $\omega = (\omega_j)_{j \in \ZZ^d}$ of (almost surely) positive, bounded, independent and  identically distributed random variables with distribution measure $\mu$, 
as well as a compactly supported, measurable function $u\colon \RR^d \to \RR$. 
The \emph{random breather potential} is the function 
\[
V_\omega(x):= \sum_{j\in \ZZ^d} u\Big(\frac{x-j}{\omega_j}\Big), 
\] 
while the family $(H_\omega)_\omega$ with $H_{\omega}:= -\Delta +V_\omega$ is called \emph{random breather model}.
\par
Random breather potentials have been introduced in \cite{CombesHM-96}, and studied in \cite{CombesHN-01} and
\cite{KirschV-10}. However, all these papers assumed unnatural regularity conditions, excluding the most basic and standard
type of single site potential, where $u$ equals the characteristic function of a ball or a cube.  
This was not a coincidence but a consequence of the linearization method used in the proofs. Our proof does not rely on linearization,
but merely on monotonicity.
While we have results for a broad class of random breather models, we restrict ourselves in this note for the purpose of clarity to
the two mentioned cases, i.e.
\begin{subequations}\label{eq:standardRBP}
\begin{align}
u&=\chi_{B_{1}},        \qquad \text{ thus }     & V_\omega(x) = \sum_{j \in \ZZ^d} \chi_{B_{\omega_j}}(x-j),  \label{eq:ballRBP}
\\
u&=\chi_{\Lambda_2},    \qquad \text{ thus }     & V_\omega(x) = \sum_{j \in \ZZ^d} \chi_{\Lambda_{2 \omega_j}}(x-j)   \label{eq:cubeRBP}
\end{align}
\end{subequations}
In fact, since our proofs are based on the analysis of level sets of random potentials, 
they work also for other types of stochastic fields with non-linear, monotone randomness,
not just for random breather potentials. 
Specifically, the function $\omega_j \mapsto \langle \phi, V_\omega\phi\rangle$ merely needs to be polynomially increasing.
\par
Note also that the random potential is uniformly bounded and non-nega\-tive, and thus the operator $H_\omega$ is self-adjoint.
\begin{theorem}[Wegner estimate for the random breather model]\label{thm:Wegner}
Let $H_\omega$ be as in \eqref{eq:standardRBP}. Assume that $\mu$ has a 
bounded density $\nu$ supported in $[\omega_{-}, \omega_{+}]$ with $0 \leq \omega_{-} < \omega_{+} < 1/2$. 
Fix $E_0 \in \RR$. Then there are
$C=C(d,E_0)$ and $\epsilon_{\max}=\epsilon_{\max}(d,E_0  ,\omega_{+}) \in (0,\infty)$
such that for all $\epsilon \in (0,\epsilon_{\max}]$ and $E \geq 0$ with
$[E-\epsilon, E+\epsilon] \subset (- \infty , E_0]$,
we have 
\begin{equation*}
\EE \left[ \mathrm{Tr} \left[ \chi_{[E- \epsilon, E + \epsilon]}(H_{\omega,L}) \right] \right]
\leq 
C  
\lVert \nu \rVert_\infty
\epsilon^{[K_0(2+{\lvert E_0 + 1 \rvert}^{1/2})]^{-1}} 
\left\lvert\ln \epsilon \right\rvert^d L^d.
\end{equation*}
The constant $\epsilon_{\max}$ can be chosen as
\begin{equation*}
\epsilon_{\max} =
\frac{1}{4} 
\left( \frac{1/2 - \omega_{+}}{2} \right)^{K_0(2+{\lvert E_0+1 \rvert}^{1/2})},
\end{equation*}
where $K_0$ is the constant from Theorem~\ref{thm:sfuc}.
\end{theorem}
Here $\EE$ denotes the expectation w.r.t.~the random variables $\omega_j, j\in \ZZ^d$, and 
$H_{\omega,L}$ the restriction of $H_{\omega}$ to the cube $\Lambda_L$ with Dirichlet boundary conditions.
Theorem \ref{thm:Wegner} implies local H\"older continuity  of the integrated density of states (IDS) and
is sufficient for the multiscale-analysis proof of spectral localization. This will be elaborated in detail elsewhere.
\begin{remark}
 The proof of Theorem~\ref{thm:sfuc} relies on Carleman estimates with and without boundary term, 
see e.g.\ \cite{LebeauR-95} and \cite{EscauriazaV-03,BourgainK-05}, on  
interpolation inequalities and an auxiliary Cauchy problem in $d+1$ dimensions
as discussed in \cite{LebeauR-95,LebeauZ-98,JerisonL-99}, and finally on geometric covering arguments developed in the theory of random Schr\"odinger operators, e.g.\ \cite{RojasMolinaV-13}.
\par
The proof of Theorem~\ref{thm:Wegner} relies on the method outlined in 
\cite{HundertmarkKNSV-06} and \cite{RojasMolinaV-13}. It can be traced back to Wegner's original work \cite{Wegner-81}.
Additional steps are necessary, since the breather model has a non-linear dependence on the random variables, unlike the well-studied Anderson model.
We also do not have the differentiability of the map  $\omega_j \mapsto \langle \phi, H_\omega\phi\rangle$  in the usual sense. 
Thus, for instance the proofs of \cite{CombesHK-07,Klein-13} do not apply.
However, the strategy of \cite{HundertmarkKNSV-06,RojasMolinaV-13} is quite versatile and can be adapted to our setting.
The key idea is not to rely on differentiability of quadratic forms but rather directly on the Courant-Hilbert  variational principle for eigenvalues.
\end{remark}
In particular, the following lemma is crucial for the proof of Theorem~\ref{thm:Wegner}. It relies on the quantitative version of the uncertainty principle from Theorem~\ref{thm:sfuc}.
Denote the eigenvalues of $H_{\omega, L}$ by $\left\{ E_n(\omega) \right\}_{n\in \NN^\ast}$, enumerated increasingly and counting multiplicities.
For $\delta \in \RR$ we define $\omega + \delta$ by $
\left( \omega + \delta \right)_j := \omega_j + \delta$ for all $j \in \ZZ^d$.
\begin{lemma}
	Let $H_{\omega,L}$ be as above and assume that $\omega \in [\omega_{-}, \omega_{+}]^{\ZZ^d}$, $\delta \leq 1/2-\omega_+$. 
	Then, for all $n \in \NN^\ast$ with $E_n(\omega) \in (- \infty , E_0]$ we have
	\begin{equation*}
	E_n(\omega + \delta) \geq E_n(\omega) + \left(\frac{\delta}{2} \right)^{\bigl[K_0 \bigl(2+ \lvert E_0+1 
	\rvert^{1/2} \bigr)\bigr]},
	\end{equation*}
	where $K_0$ is the constant from Theorem~\ref{thm:sfuc}.
	\end{lemma}
Thus, we obtain a lifting estimate on the eigenvalues which is independent on the length scale.
Details of the proof of Theorem~\ref{thm:Wegner} can be found in \cite{TaeuferV}.
\begin{remark}[Challenges due to non-linearity]\label{r-nonlinearity}
The challenges are best understood by comparing the breather model with the 
alloy-type potential $V_\omega (x) =\sum_{j\in \ZZ^d} \omega_j u(x-j)$ (for simplicity consider $u=\chi_{B(r)}$). 
The latter  depends in a linear way on the random coupling constants constituting the configuration 
$\omega = (\omega_j)_j$. In particular, the derivatives of eigenvalues $E_n(\omega)$ (of finite box restrictions on $-\Delta +V_\omega)$
w.r.t.~each $\omega_j$ are easily calculated 
via the Hellman-Feynman formula.
\begin{figure}[t]\centering
\includegraphics[width=12cm]{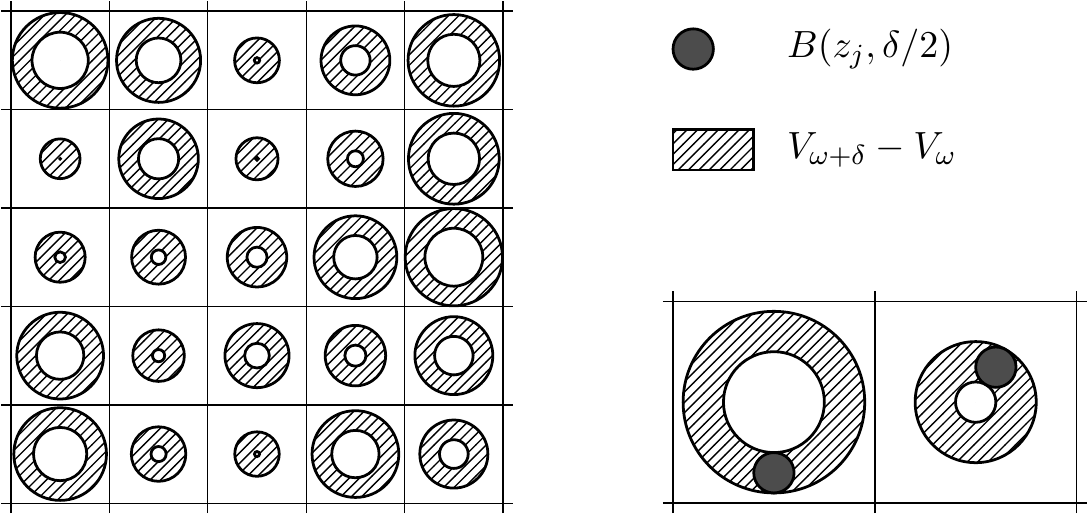}

\caption{Illustration of the support of the increments $V_{\omega + \delta} - V_\omega$ (left) 
and the choice of the balls $B(z_j, \delta/2)$ (right).
(Illustration des support des incr\'ements $V_{\omega + \delta} - V_\omega$ et du choix des boules $B(z_j, \delta/2)$)\label{fig:annuli}}
\end{figure}
In contrast, for the breather model the derivatives $\frac{\partial}{\partial \omega_j}E_n(\omega)$ are only defined in distributional sense.
Thus, one is lead to implement eigenvalue perturbation theory using increments $V_{\omega+\delta}-V_{\omega}$, with positive $\delta$.
Note that in the case of the alloy-type model, for any fixed $\delta$, the increment $V_{\omega+\delta}-V_{\omega}$ is independent of the configuration 
$\omega$ and a $\ZZ^d$-periodic function. Therefore, it is not needed to know the explicit dependence of $C_{\rm sfuc}$ on $\delta$. For the breather model this is not the case. In particular, $V_{\omega+\delta}-V_{\omega}$ is a non-periodic function and its support depends both on $\delta$ and $\omega$. Specifically, it is a union of annuli of width $\delta$ and $\omega$-dependent radii, cf.\ Fig.~\ref{fig:annuli}.
Technically, one has to estimate the mass of the square of an eigenfunction in this support set as a function of $\omega$ and $\delta$.
For the application of Theorem 2.1 one has to chose balls $B (z_j, \delta / 2)$  lying inside the annuli, see Fig.~\ref{fig:annuli}.
To obtain H\"older continuity of the IDS one has to control the behaviour of $C_{\rm sfuc}$ as $\delta \searrow 0$.
\end{remark}
\def\cprime{$'$}\def\polhk#1{\setbox0=\hbox{#1}{\ooalign{\hidewidth
  \lower1.5ex\hbox{`}\hidewidth\crcr\unhbox0}}}


\begin{thebibliography}{10}
%
\bibitem{BourgainK-05}
J.~Bourgain and C.~E. Kenig.
\newblock On localization in the continuous {A}nderson-{B}ernoulli model in
  higher dimension.
\newblock {\em Invent. Math.}, 161(2):389--426, 2005.

\bibitem{BourgainK-13}
J.~Bourgain and A.~Klein.
\newblock Bounds on the density of states for {S}chr\"odinger operators.
\newblock {\em Invent. Math.}, 194(1):41--72, 2013.

\bibitem{CombesHK-07}
J.-M. Combes, P.~D. Hislop, and F.~Klopp.
\newblock An optimal {Wegner} estimate and its application to the global
  continuity of the integrated density of states for random {Schr\"odinger}
  operators.
\newblock {\em Duke Math. J.}, 140(3):469--498, 2007.

\bibitem{CombesHM-96}
J.-M. Combes, P.~D. Hislop, and E.~Mourre.
\newblock Spectral averaging, perturbation of singular spectra, and
  localization.
\newblock {\em Trans. Amer. Math. Soc.}, 348(12):4883--4894, 1996.

\bibitem{CombesHN-01}
J.-M. Combes, P.~D. Hislop, and S.~Nakamura.
\newblock The ${L}^p$-theory of the spectral shift function, the {Wegner}
  estimate, and the integrated density of states for some random
  {Schr\"odinger} operators.
\newblock {\em Commun. Math. Phys.}, 70(218):113--130, 2001.

\bibitem{EscauriazaV-03}
L.~Escauriaza and S.~Vessella.
\newblock Optimal three cylinder inequalities for solutions to parabolic
  equations with {L}ipschitz leading coefficients.
\newblock In {\em Inverse problems: theory and applications ({C}ortona/{P}isa,
  2002)}, volume 333 of {\em Contemp. Math.}, pages 79--87. Amer. Math. Soc.,
  Providence, RI, 2003.

\bibitem{HundertmarkKNSV-06}
D.~Hundertmark, R.~Killip, S.~Nakamura, P.~Stollmann, and I.~Veseli\'c.
\newblock Bounds on the spectral shift function and the density of states.
\newblock {\em Commun. Math. Phys.}, 262(2):489--503, 2006.

\bibitem{JerisonL-99}
D.~Jerison and G.~Lebeau.
\newblock Nodal sets of sums of eigenfunctions.
\newblock In {\em Harmonic analysis and partial differential equations
  ({C}hicago, {IL}, 1996)}, Chicago Lectures in Math., pages 223--239. Univ.
  Chicago Press, Chicago, IL, 1999.

\bibitem{KirschV-10}
W.~Kirsch and I.~Veseli\'c.
\newblock Lifshitz tails for a class of Schr\"odinger operators with random
  breather-type potential.
\newblock {\em Lett. Math. Phys.}, 94(1):27--39, 2010.

\bibitem{Klein-13}
A.~Klein.
\newblock Unique continuation principle for spectral projections of
  Schr{\"o}dinger operators and optimal Wegner estimates for non-ergodic random
  Schr{\"o}dinger operators.
\newblock {\em Commun. Math. Phys.}, 323(3):1229--1246, 2013.

\bibitem{RousseauL-12}
{J.~Le}~Rousseau and G.~Lebeau.
\newblock On {C}arleman estimates for elliptic and parabolic operators. {A}pplications to unique continuation and control of parabolic equations.
\newblock {\em ESAIM Contr. Optim. Ca.},
  18(3):712--747, 2012.
	
\bibitem{LebeauR-95}
G.~Lebeau and L.~Robbiano.
\newblock Contr{\^o}le exact de l{\'e}quation de la chaleur.
\newblock {\em Commun. Part. Diff. Eq.},
  20(1-2):335--356, 1995.

\bibitem{LebeauZ-98}
G.~{Lebeau} and E.~{Zuazua}.
\newblock Null-controllability of a system of linear thermoelasticity.
\newblock {\em Arch. Ration. Mech. An.}, 141(4):297--329,
  1998.

\bibitem{NakicTTV-prep}
I.~Naki\'c, M.~T\"aufer, M.~Tautenhahn, and I.~Veseli\'c.
\newblock In preparation.
  
\bibitem{RojasMolinaV-13}
C.~Rojas-Molina and I.~Veseli\'c.
\newblock Scale-free unique continuation estimates and applications to random
  Schr\"odinger operators.
\newblock {\em Commun. Math. Phys.}, 320(1):245--274, 2013.

\bibitem{TaeuferV}
M.~T\"aufer and I.~Veseli\'c
\newblock  Conditional Wegner estimate for the standard random breather potential.
\newblock {\em J. Stat. Phys.}. DOI/10.1007/s10955-015-1358-y, 2015.

\bibitem{Wegner-81}
F.~Wegner.
\newblock Bounds on the {DOS} in disordered systems.
\newblock {\em Z. Phys. B}, 44:9--15, 1981.

\end{thebibliography}
\end{document}